\RequirePackage{fix-cm}
\documentclass[a4paper,english,10pt,oneside]{article}\def\zibreport{1}
\def\arxiv{1}
\usepackage{graphicx}

\makeatletter
\let\cl@chapter\undefined
\makeatletter
\usepackage{mathptmx}      

\usepackage{ifthen}
\usepackage{multicol}
\usepackage{footmisc}
\usepackage{mathdots}
\usepackage{amsmath,amsfonts,amssymb,mathrsfs}
\usepackage{booktabs}
\usepackage{xspace}
\usepackage{lscape}
\usepackage{multirow}
\usepackage{comment}
\usepackage{hyperref}
\usepackage{cleveref}

\usepackage{url}

\usepackage{orcidlink}

\usepackage{verbatim}
\usepackage{algorithm}
\usepackage{algpseudocodex}
\usepackage{todonotes}

\newcommand{\myorcidlink}[1]{\,\href{https://orcid.org/#1}{\raisebox{-0.45ex}{\includegraphics[width=1.8ex]{orcid}}}}
\newcommand{\myurl}[1]{\textsf{\footnotesize \url{#1}}\xspace}

\newcommand{\allcaps}[1]{\protect\scalebox{0.93}{#1}}

\newcommand{\name}[1]{\mbox{#1}\xspace}
\newcommand{\nameCaps}[1]{\name{\allcaps{#1}}}

\newcommand{\remix}{\allcaps{REMix}\xspace}

\newcommand{\UNSEEN}{\nameCaps{UNSEEN}}

\newcommand{\CPLEX}{\nameCaps{CPLEX}}

\newcommand{\GUROBI}{\name{Gurobi}}

\newcommand{\PDLP}{\nameCaps{PDLP}}
\newcommand{\PDHG}{\nameCaps{PDHG}}

\newcommand{\COPT}{\nameCaps{COPT}}

\newcommand{\RAM}{\nameCaps{{RAM}}}
\newcommand{\GB}{\nameCaps{{GB}}}
\newcommand{\ARM}{\nameCaps{{ARM}}}
\newcommand{\AMD}{\nameCaps{{AMD}}}
\newcommand{\EPYC}{\nameCaps{{EPYC}}}
\newcommand{\NVIDIA}{\nameCaps{{NVIDIA}}}

\newcommand{\GPUs}{\nameCaps{GPUs}}
\newcommand{\GPU}{\nameCaps{GPU}}
\newcommand{\CPUs}{\nameCaps{CPUs}}
\newcommand{\CPU}{\nameCaps{CPU}}
\newcommand{\MIPLIB}{\nameCaps{MIPLIB}}

\newcommand{\LP}{\nameCaps{LP}}
\newcommand{\LPs}{\nameCaps{LPs}}
\newcommand{\MIP}{\nameCaps{MIP}}
\newcommand{\FP}{\nameCaps{FP}}

\newcommand{\MIPs}{\nameCaps{MIPs}}

\newcommand{\ESOM}{\nameCaps{ESOM}}
\newcommand{\ESOMs}{\nameCaps{ESOMs}}

\newcommand{\IPM}{\nameCaps{IPM}}
\newcommand{\IPMs}{\nameCaps{IPMs}}

\newcommand{\NP}{\nameCaps{NP}}

\newcommand{\FOMs}{\nameCaps{FOMs}}

\newcommand{\var}{\ensuremath{x}\xspace}

\newcommand{\primalvector}{\var}

\newcommand{\objectiveLinearVector}{\ensuremath{c}\xspace}

\newcommand{\minPrimalIn}[1]{\min_{\primalvector \in {#1}}}

\newcommand{\R}{\mathbb{R}}
\newcommand{\Rtimes}[2]{\R^{{#1}\times{#2}}}

\newcommand{\Rmn}{\Rtimes{m}{n}}

\newcommand{\Z}{\mathbb{Z}}

\definecolor{seagreen}{rgb}{0.18,0.74,0.56}
\definecolor{darkgreen}{rgb}{0.0,0.45,0.00}
\definecolor{navyblue}{rgb}{0.0,0.0,0.5}
\definecolor{steelblue}{rgb}{0.27,0.51,0.71}
\definecolor{siennabrown}{rgb}{0.63,0.32,0.18}
\definecolor{firebrickred}{rgb}{0.69,0.13,0.13}
\definecolor{gray75}{rgb}{0.75,0.75,0.75}
\definecolor{orange}{rgb}{.843,0.671,0.078}
\definecolor{gold}{rgb}{1.0,0.84,0.0}
\definecolor{scipyellow}{HTML}{FFFFD6}
\definecolor{soplexred}{HTML}{FFD8D8}
\definecolor{zimplgreen}{HTML}{D8FFD8}
\definecolor{ugblue}{HTML}{CFEFFF}
\definecolor{gcgorange}{HTML}{FFDAB9}

\definecolor{mDarkBrown}{HTML}{604c38}
\definecolor{mDarkTeal}{HTML}{23373b}
\definecolor{mLightBrown}{HTML}{EB811B}
\definecolor{mLightblue}{HTML}{14B03D}

\definecolor{c0}{HTML}{000060}
\definecolor{c1}{HTML}{0000FF}
\definecolor{c2}{HTML}{36648B}
\definecolor{c3}{HTML}{4682B4}
\definecolor{c4}{HTML}{5CACEE}
\definecolor{c5}{HTML}{FF0000}
\definecolor{c6}{HTML}{008888}
\definecolor{c7}{HTML}{00DD99}
\definecolor{c8}{HTML}{527B10}
\definecolor{c9}{HTML}{7BC618}
\definecolor{c10}{HTML}{8DD8F8}
\definecolor{background}{HTML}{FFFFFF}

\definecolor{cp1}{HTML}{F2AF29}
\definecolor{cp2}{HTML}{05668D}
\definecolor{cp3}{HTML}{02C39A}
\definecolor{cp4}{HTML}{4F345A}
\definecolor{cp5}{HTML}{F2B5D4}
\definecolor{cp6}{HTML}{DA2C38}
\DeclareMathAlphabet{\mathcal}{OMS}{cmsy}{m}{n}

\newcommand{\citet}{\cite}
\newcommand{\citep}{\cite}

\ifthenelse{\zibreport = 1}{}{
    \smartqed  
}

\newcommand{\TheTitle}{Fix-and-Propagate Heuristics Using Low-Precision First-Order LP Solutions for Large-Scale Mixed-Integer Linear Optimization} 
\newcommand{\TheAuthors}{N.--C. Kempke, T. Koch}

\hypersetup{%
    pdftitle={\TheTitle},
    pdfauthor={\TheAuthors}
}

\ifthenelse{\zibreport = 1}{
    \usepackage{./zibtitlepage}
}{
    \journalname{Mathematical Programming Computation}
}

\begin{document}

\ifthenelse{\zibreport = 1}{
    \ZTPTitle{\TheTitle}
    \title{\TheTitle}

    \ZTPAuthor{
        \ZTPHasOrcid{Nils--Christian Kempke}{0000-0003-4492-9818},
        \ZTPHasOrcid{Thorsten Koch}{0000-0002-1967-0077}
    }
    \author{
        \ZTPHasOrcid{Nils--Christian Kempke}{0000-0003-4492-9818},\and\
        \ZTPHasOrcid{Thorsten Koch}{0000-0002-1967-0077}
    }

    \ZTPInfo{Preprint}
    \ZTPNumber{25-03}
    \ZTPMonth{March}
    \ZTPYear{2026}

    \date{\normalsize March 3, 2026}
    \ifthenelse{\arxiv = 0}{
        \zibtitlepage
    }{}
    \maketitle

}{
    \title{\TheTitle
    }

    \titlerunning{Fix-And-Propagate via Low-Precision \FOMs}        

    \author{Nils-Christian Kempke\orcidlink{0000-0003-4492-9818} \and
        Thorsten Koch\orcidlink{0000-0002-1967-0077}
    }


    \institute{Nils-Christian Kempke (corresponding author) \and Thorsten Koch \at
                Zuse Institute Berlin, Applied Algorithmic Intelligence Methods, Takustre{\ss}e~7, 14195~Berlin, Germany \\
                \email{\{kempke,koch\}@zib.de}
            \and
            Thorsten Koch \at Technical University Berlin, Chair of Software and Algorithms for Discrete Optimization, Stra{\ss}e des 17.~Juni~136, 10623~Berlin, Germany
    }

    \date{Received: date / Accepted: date}
}

\maketitle
\begin{abstract}
We investigate the use of low-precision first-order methods (\FOMs) within a fix-and-propagate (\FP) framework for solving mixed-integer programming problems (\MIPs).
We employ \GPU-accelerated \PDLP, a variant of the Primal-Dual Hybrid Gradient (\PDHG) method specialized to \LP problems, to solve the \LP-relaxation of our \MIPs to low accuracy.
This solution is used to motivate fixings within our \FP framework.
We evaluate the performance of our heuristic on \MIPLIB 2017, demonstrating that low-accuracy \LP solutions do not lead to a loss in the quality of the \FP heuristic solutions.
Further, we use our \FP framework to produce high-accuracy solutions for large-scale (up to 243 million nonzeros and 8 million decision variables) unit commitment-based dispatch and expansion planning problems created with the modeling framework \remix.
For the largest problems, we can generate solutions with a primal-dual gap of under 2\% in less than 4 hours, whereas state-of-the-art commercial solvers cannot produce feasible solutions within 2 days of runtime.
\ifthenelse{\zibreport = 1}{}{
}
\end{abstract}

\section{Introduction}

Linear programming problems (\LPs) optimize a linear objective subject to linear constraints:
\begin{equation} \label{eq:LP}
  \minPrimalIn{\R^n} \{\objectiveLinearVector^T \var : A \var \leq b, l \leq \var \leq u\},
\end{equation}
where $\var$ are decision variables, $c\in\R^n$ objective, $A\in\Rmn$ the constraint matrix, $b\in\R^m$ the right-hand sides, and $l,u \in \overline{\R} := \R \cup \{\pm \infty\}$ the variable bounds.
\LPs are typically solved with Simplex \citep{DantzigWolfeSimplex_1955} or Interior-Point Methods (\IPMs) \citep{wright1997_primal_dual_ipm}.
Recently, First-Order Methods (\FOMs), using gradient information and matrix-vector products, gained attention for large-scale \LPs due to their \GPU suitability \cite{ApplegatePDLP,LoRADS_Han2024,ABIP_ADMM1}.

A mixed-integer program (\MIP), ex\-tends \cref{eq:LP} with integrality constraints on a subset of variables, the \textit{integer variables} or short \textit{integers}: $\var_i \in \Z\ \forall\ i\in\mathcal{I}\subset \{1,\dots,n\}.$
\MIPs are \NP-hard and the most established solution methods are branch-and-bound approaches relying on \IPMs and Simplex methods.
Additionally, general purpose \MIP solvers heavily rely on primal heu\-ris\-tics \citet{PhDAchterberg} generating feasible solutions quickly.
Fix-and-propagate (\FP) heuristics \citep{Achterberg2013} iteratively fix integers and apply domain propagation to reduce the domain of other variables.
If the \FP heuristic finds a feasible integer assignment, all integers are fixed and the resulting \LP is solved generating a \MIP feasible solution or aborting the heuristic if the \LP is infeasible.

Among \FOMs, \PDLP has recently been integrated into the Cardinal Optimizer (\COPT), Fico Xpress, NVIDIA cuOpt, and Gurobi, enabling \GPU accelerated \LP solves and providing \PDLP as an alternative method for the root \LP relaxation of \MIPs.
Apart from this, little \GPU integration has been done into commercial/academic software and this raises the question: can \FOMs accelerate heuristics within \MIP frameworks while preserving solution quality?
Prior work includes Scylla \cite{Scylla_Mexi2023}, embedding CPU-based \PDLP in context of the feasibility pump \citep{FeasPump_Fischetti2005} where already in \citet{FeasPump_Fischetti2005}, the possibility of increasing the speed of the heuristic via low-quality \LP solves is entertained.
During the revision process, cuOpt added a \GPU-based feasibility pump using \PDLP in combination with propagation and local search \cite{cuOptFeasPump}, as the first commercial standalone \GPU-accelerated primal heuristic.
Additionally, \COPT and cuOpt have both added a \GPU accelerated \IPM.

\paragraph{Contribution}

We make three contributions:
\begin{itemize}
    \item A \FP framework exploiting \LP solutions of varying accuracy to generate high-quality \MIP solutions.
    \item Performance evaluation on \MIPLIB~2017 \citep{MIPLIB2017}, showing that \LP solution accuracy barely affects performance while allowing for speed-ups one large models.
    \item Application to large-scale energy system optimization models (\ESOMs), closing gaps $< 2\%$ within 4 hours, showing speedups up to a factor of 3 over the \IPM.
\end{itemize}

We show with our \MIPLIB experiments, that \FOMs can readily be embedded into fix-and-propagate type heuristics without sacrificing solution quality.
With our \ESOM experiments we highlight a class of models where this embedding can be beneficial.

\paragraph{Outline}

\Cref{sec:pdlp_fap}, details the heuristic framework \LP integration strategies.
\Cref{sec:models} introduces our \ESOMs.
\Cref{sec:numerical_results} presents our \MIPLIB and \ESOM results.
\section{Fix-and-Propagate Using Low-Precision LP Solutions}\label{sec:pdlp_fap}

Our implementation builds on the depth-first-serach (\textit{dfs} and \textit{dfsrep}) mode of the fix-propagate-repair framework from \cite{Salvagnin2024_AFixPropagateRepairHeuristicsForMIP} where the main \FP scheme is described in ``Fig. 1 Fix-propagate-repair scheme''.
Before running, the framework uses a \textit{variable strategy} sorting the integers and binaries.
It then starts the fix-propagate-repair process using a \textit{value strategy} to determine fixing values, fixing variables in the order computed before.
We used all propagators from the original code: clique, implication, and linear constraint propagation.
We extended the original code in three ways:
\begin{itemize}
    \item We added \LP-based variable strategies (extending Table 1 in \cite{Salvagnin2024_AFixPropagateRepairHeuristicsForMIP}): \textit{Frac}, \textit{Duals}, \textit{RedCost}, and \textit{Type} (and several hybrid strategies).
    \item We added a variable branching strategy better suited for integers, replacing the call in line 14 of Fig.~1 in \cite{Salvagnin2024_AFixPropagateRepairHeuristicsForMIP}.
    \item We enabled the use of \PDLP (or any other \LP solver) for solving \LP relaxations within the framework with different tolerances.
\end{itemize}

For any \LP-based variable or value strategy, the framework first solves an \emph{initial \LP} using Simplex, \IPM, or \PDLP.
The relative tolerance of the initial solve for \PDLP and \IPM can be adjusted, for Simplex we use solver's default tolerances.
The \LP solution guides variable selection and/or fixing throughout the \FP branching procedure.
Periodic re-solving to improve the \LP solution information is possible in principle but not used here, as it incurs substantial computational cost.
After finding a partial assignment of integers, the \emph{final \LP} is solved with an \IPM to produce a feasible solution.

\paragraph{\LP-based variable strategies}

\textit{Frac} sorts variables by fractionality, the distance of the \LP value to the nearest integer: $ \text{frac}_i = \min(x_i - \lfloor x_i \rfloor, \lceil x_i \rceil - x_i) < 1$.
This strategy is commonly used in diving or \FP heuristics \cite{MasterThesisTimo_2006}.

\textit{RedCost} uses a variable's reduced costs to determine its ordering.
The reduced cost indicates how much the objective would change if the variable moved from its current value.
Variables with high absolute reduced cost are less likely to deviate in an optimal solution (assuming the optimum is close to the \LP solution value).
Greedily fixing these variables first aims to stay close to the objective value.

\textit{Dual} combines dual constraint values with reduced costs.
Intuitively, constraints with large absolute duals often have a strong influence on the objective and tend to be tight in optimal solutions.
We sort constraints by absolute dual value, iterate through them, and extract integer variables not yet considered, ordering these by reduced cost.
This heuristic prioritizes variables likely to impact the objective associated with tight constraints and follows an intuition rather than a formally proven property.

\textit{Type} uses the pre-existing variable strategy \textit{type} sorting variables stably, putting binaries first and integers last.

\paragraph{LP-based fixing}

In our experiments we pair \textit{Frac}, \textit{Redcosts}, \textit{Dual}, and \textit{Type} with an \LP-based fixing strategy.
Given an integer variable $x_i$ with \LP value $x_{\text{\LP},i}$, we compute its fractional distance to the lower integer: $d_i=x_{\text{\LP},i} - \lfloor x_{\text{\LP},i} \rfloor$.
We then sample $d \sim U(0,1)$ and fix $x_i$ to its lower bound if $d > d_i$, otherwise to its upper bound.
This choice of fixing is meant to stay close to the \LP solution while also allowing for some diversification of the search space and avoiding deterministic patterns when running the heuristic more than once or with similar \LP solutions.
Fractionality is treated as a ``confidence'' of the variable being a particular integer or not and for values close to $0.5$ the strategy becomes random.
Strictly rounding to the closest integer usually leads to infeasibility and, know as \textit{simple rounding}, is also discussed in \cite{MasterThesisTimo_2006}.

\paragraph{Integer-aware branching strategy}

Our branching strategy, shown in \cref{alg:branching}, generates two or three \textit{branching bound changes} depending on the domain of the integer variable $x_i$ and its fixing value $a_i \in \mathbb{Z}$.
Each branching bound change consists of the new lower and upper bound that should be enforced on a variable and the variable's index.
If the fixing value matches a bound (line 6), we create two branching bound changes (line 8 or 10), fixing the variable to either bound, preferring the one equal to the fixing value.
Otherwise, three nodes are generated: one fixing the variable to the suggested value, and two restricting its domain above and below it.
When generating three nodes (line 3 or 5) we prefer fixing to the branching value, and then the restriction in objective direction.
These restrictions do not fix the variable but ensure it is soon branched on again.
\LP-based fixing often suggests values inside the domain, helping explore regions around \LP solutions.  
This change was especially important for our \ESOMs as shown in detail in \cref{sec:appendix:esoms_res}.

\begin{algorithm}[h!]
\caption{Branching strategy}\label{alg:branching}
\begin{algorithmic}[1]
\Require Integer variable $x_i$ with bounds $l_i \leq x_i \leq u_i $, branching value $a_i \in \mathbb{Z}$ with $l_i \leq a_i \leq u_i$
\Ensure $({l}^1, {u}^1, i),\dots,({l}^k, {u}^k, i)$ $k$ branching bound changes in \textbf{increasing} priority
\If{$l_i < a_i < u_i$}
    \If{$c_i > 0$}
        \State \Return $[(a_j + 1, u_j, i), (l_i, a_i -1, i), (a_j, a_j, i)]$
    \Else
        \State \Return $[(l_i, a_i -1, i), (a_j + 1, u_j, i), (a_j, a_j, i)]$
    \EndIf
\Else
    \If{$a_i = l_i$}
        \State \Return $[({u}_j, {u}_j), ({l}_j, {l}_j)]$ 
    \Else
        \State \Return $[({l}_j, {l}_j), ({u}_j, {u}_j)]$ 
    \EndIf
\EndIf
\end{algorithmic}
\end{algorithm}

\paragraph{Tiebreaking and hybrid strategies}

When two or more variables are ranked equally by a given variable strategy, a tiebreaking strategy can chose the potentially ``better'' variable for fixing next.
In \cref{sec:numerical_miplib} we considered \textit{Frac} with ties broken by \textit{Dual}, \textit{Frac} with ties broken by \textit{RedCost}, and vice versa. 
These hybrid approaches do not fundamentally change the underlying heuristic, but provide additional flexibility in variable selection and help avoid arbitrary decisions when multiple variables have similar \LP characteristics.
We note that the original code did not feature any tiebreaking and we run all strategies in \cite{Salvagnin2024_AFixPropagateRepairHeuristicsForMIP} Table~3 without additional tiebreaking.
\section{Energy system optimization models}\label{sec:models}

The \ESOMs used in the experiments in \cref{sec:numerical_results} were developed in the \UNSEEN project and are publicly available \cite{ZenodoUNSEEN}; all sets are used except miso\_244k.
Based on the \textsc{remix} modeling framework \citep{Wetzel2024_REMIX}, the instances represent the German power system and neighboring countries and optimize the power sector for 2030.

Existing conventional capacities follow the current fleet, while a CO$_2$ price incentivizes renewable expansion. Renewable generation (wind and solar), lithium-ion batteries, and pumped-storage hydropower are modeled with continuous variables; discrete technologies and transmission expansions use integer variables.
Natural-gas CCGT is the only expandable conventional technology; lignite, coal, and OCGT capacities are fixed.
Unit commitment includes minimum up-/down-times and partial-load constraints for CCGT, coal, and lignite.
The transmission network comprises HVAC and HVDC corridors.
AC lines use a DC load-flow (B--$\theta$) formulation, and HVDC links are modeled as capacity-constrained transport arcs.
Grid topology is fixed and expandable only in capacity.
Investments are static, chosen once for the modeled year; commissioning times and staged investments are not considered.

The resulting large-scale \MIPs span 8{,}760 hourly time steps.
Model sizes and instance counts are reported in \cref{tab:unseen_sizes}.
\begin{table}[ht!]
\caption{Sizes and number of \UNSEEN \ESOMs.}
\label{tab:unseen_sizes}
\begin{tabular}{l | r | r | r | r | r} 
\toprule
& Variables & Integers & Constraints & Non-zeros & \#instances\\ 
\midrule
miso\_2M & 880,606 & 35,052 & 928,790 & 2,829,975 & 100 \\
miso\_4M & 1,103,903 & 157,692 & 1,279,124 & 4,590,593 & 96 \\ 
miso\_47M & 13,597,290 & 972,477 & 16,313,022 & 47,310,667 & 100 \\
miso\_82M & 24,356,961 & 867,465 & 24,698,935 & 82,773,747 & 20 \\
remix\_105M & 28,455,474 & 3,364,138 & 33,194,746 & 105,960,838 & 100 \\
remix\_243M & 66,714,360 & 7,726,975 & 77,498,047 & 243,495,063 & 100 \\
\bottomrule
\end{tabular}
\end{table}
These undecomposed formulations capture strong seasonal effects in renewable-dominated systems (see \cite{CaoTempZooming} for comparisons with rolling-horizon and temporal-zooming approaches).

\section{Numerical results}\label{sec:numerical_results}

\paragraph{Hardware and software} The \MIPLIB experiments in \cref{sec:numerical_miplib} ran on four identical machines, each with four \NVIDIA A100 \GPUs (80 \GB), 128 \GB \RAM, and two \AMD \EPYC~7513 \CPUs (32 cores each).
One \CPU socket was allocated per run, with up to two concurrent runs per machine, all runs being allocated 32 threads.
The \ESOM experiments in \cref{sec:numerical_results_esom_pdlp_vs_ipm,sec:numerical_results_ESOM} used a single \NVIDIA GH200 Grace--Hopper system with a 72-core \ARM Neoverse-V2 \CPU, 480 \GB \CPU memory, and one \NVIDIA H200 \GPU (96 \GB).
All runs used 32 threads (affecting only the \IPM); two parallel runs were allowed, but never simultaneously using the \GPU.

Instances were presolved with \GUROBI~11, which is essential for tractability; presolve time is included in total runtimes.
Only presolved solutions are reported, as postsolve is not exported by \GUROBI.
Although the framework also supports \CPLEX for pre- and postsolve, but when starting this paper \CPLEX offered only limited support on Linux \ARM.
Open-source presolvers such as PaPILO were considered but performed poorly on our \ESOMs.
We use \COPT for \PDLP, \IPM, and Simplex solves using \COPT~7.2.0 on x86 and 7.2.3 on \ARM.

\paragraph{FP experiments} We conducted three experiments: (i) in \cref{sec:numerical_miplib}, each \FP strategy, including the presets from Table~3 \cite{Salvagnin2024_AFixPropagateRepairHeuristicsForMIP}, was run on the \MIPLIB~2017 benchmark \citep{MIPLIB2017}; after removing seven infeasible instances and applying five random permutations per feasible instance \citep{PerfVariabilityInMIP_Lodi2013}, this yielded 1{,}165 instances; (ii) in \cref{sec:numerical_results_esom_pdlp_vs_ipm}, we compared \IPM and \PDLP runtimes on large-scale \ESOMs; and (iii) in \cref{sec:numerical_results_ESOM}, we evaluated \LP-based strategies using \PDLP and \IPM on the \ESOMs.

The \IPM is used for the final \LP relaxation with primal, dual, and gap shiftetolerances of $10^{-8}$.
The \LP method for the initial relaxation depends on the strategy: Dual Simplex (presets \textit{zerolp} and \textit{lp}) uses \COPT defaults of $10^{-6}$ absolute primal and dual tolerances, while \PDLP uses relative tolerances.
\IPMs often yield absolute primal-dual feasibility even for larger gap tolerances, whereas \FOMs provide only relative feasibility.
A solution is \textit{feasible} if it satisfies relative primal/dual tolerances of $10^{-6}$.

\paragraph{Other} Unless stated otherwise, results are reported using the shifted geometric mean with a shift of one second or percent reducing sensitivity to outliers.
Gaps and solution times are reported only for instances with feasible solutions; unsolved instances typically have shorter runtimes as the final \LP is not executed.
Gaps are computed relative to known optima for \MIPLIB and to the \LP-relaxation value for \ESOMs.
Results often aggregate multiple configurations by taking the best outcome per instance, reducing the impact of performance variability \cite{PerfVariabilityInMIP_Lodi2013}.
Full tables are provided in the appendix, detailed logs in \cite{ZenodoFPLogs}.
Finally, as NVIDIA released cuOpt during revision, limited results on cuOpt for \ESOMs and \MIPLIB perturbations are reported in \cref{sec:appendix:cuopt,tab:appendix:cuopt}.

\subsection{MIPLIB experiments}\label{sec:numerical_miplib}

The \MIPLIB experiments were run with a timelimit of one hour and 20 \GB \CPU memory.
Neither limit was ever hit.

\paragraph{LP-based FP compared with LP-free variants}

We ran all presets in Table~3 \cite{Salvagnin2024_AFixPropagateRepairHeuristicsForMIP} and our additions with and without repair and with and without our newly implemented branching rules.
Full results can be found in \cref{sec:appendix:miplib_res} (\cref{tab:appendix:miplib_presets,tab:appendix:miplib_lpbased}).
In \cref{tab:miplib_lpbased_vs_lpfree} we aggregated our \LP-based variants and compare them against all presets and only \LP-free presets.
\begin{table}[ht!]
\caption{\MIPLIB results aggregated by \LP-based and \LP-free \FP strategies.}
\label{tab:miplib_lpbased_vs_lpfree}
\begin{tabular}{lccr}
\toprule
Strategies & \# found & gap (\%) & time (s) \\
\midrule
Presets            & 897 & 24.44 & 2.37 \\
Presets (\LP-free) & 849 & 33.90 & 1.35 \\
\LP based (ours)   & 840 & 16.33 & 3.74 \\
\midrule
Best                    & 922 & 15.30 & 3.66 \\
\bottomrule
\end{tabular}
\end{table}
Our additions produce fewer solutions but also result in over eight percent better gaps while also taking more time as they always rely on the initial \LP relaxation.
\LP-based methods dominate the aggregated runtime of ``Best'' and as they produce better solutions.

\paragraph{Impact of branching strategy}

\begin{table}[ht!]
	\caption{\MIPLIB results aggregated by branching strategy.}\label{tab:branching}
	\begin{tabular}{l|ccc|ccc}
		\toprule
        \multicolumn{1}{c}{} & \multicolumn{3}{c}{old branching} & \multicolumn{3}{c}{new branching} \\
		Methods  & \# found & gap (\%) & time (s) & \# found & gap (\%) & time (s) \\
		\midrule
		Presets  & 895 & 26.34 & 2.7 & 900 & 24.36 & 2.41 \\
		LP based (ours) & 832 & 18.47 & 3.65 & 834 & 16.81 & 3.77 \\
		\bottomrule
	\end{tabular}
\end{table}
To evaluate our branching strategy, we aggregated all results using the old and new branching strategies over the presets and our \LP based additions in \cref{tab:branching}.
The extended branching strategy is able to reduce the gap produced on average by about two percent, solving slightly more instances at the cost of $0.25$s higher runtimes for the presets and $0.1$s for our additions.

\paragraph{Tiebreaking/hybrid strategies}
\begin{table}[ht!]
\caption{\MIPLIB results for \LP-based and tiebreaking strategies aggregated by \LP method (tolerance $10^{-4}$).}
\label{tab:miplib_tiebreak}
\begin{tabular}{l|ccc|ccc}
\toprule
\multicolumn{1}{c}{} & \multicolumn{3}{c}{\PDLP} & \multicolumn{3}{c}{\IPM} \\
\cmidrule(lr){2-4} \cmidrule(lr){5-7}
Variable strategies & \# found & gap (\%) & time (s) & \# found & gap (\%) & time (s) \\
\midrule
\LP-based tiebreaking & 649 & 20.18 & 3.63 & 641 & 20.15 & 2.16 \\
\LP-based   & 671 & 22.86 & 3.57 & 660 & 21.92 & 2.00 \\
\midrule
Best & 693 & 19.65 & 3.67 & 699 & 19.11 & 2.19 \\
\bottomrule
\end{tabular}
\end{table}
We ran each of our tiebreaking variants using \IPM and \PDLP (tolerances $10^{-4}$ each) for the initial \LP solve and repair disabled (to decrease random influence) and compare the tiebreaking variants against their \LP counterparts.
The aggregated results are presented in \cref{tab:miplib_tiebreak}.
Overall, while the tiebreaking variants produce about 2\% smaller gaps, they also found solutions for about 20 instances less.
More notably, aggregating the tiebreaking and standard variants together, we see that they seem to complement each other, producing smaller gaps and more solutions.
Primal heuristics are strongly susceptible to variance.
By adding more variants, this randomness can be exploited to generate a more complete probing of the overall solution space which we believe it the primary benefit of the tiebreaking variants.

\paragraph{Impact of LP algorithm and precision} \Cref{tab:miplib_results_fp} shows results for our \LP-based variants using \IPM and \PDLP as \LP solver and precisions $10^{-4}$ and $10^{-6}$.  
Overall, solution quality is largely insensitive to \LP algorithm and precision for the \MIPLIB instances varying by only 0.2\% over all aggregations.
As \IPMs converge quickly and often close the final gap in few iterations for small instances (such as in \MIPLIB) reducing precision only yields minimal time savings.
\PDLP on the other hand exhibits slower final convergence but faster initial progress; reducing precision yields roughly a 1.5x speed-up.
\begin{table}[ht!]
	\caption{\MIPLIB results for \LP-based \FP variants aggregated by \LP method and tolerance.}
	\label{tab:miplib_results_fp}
	\begin{tabular}{l|ccc|ccc@{}}
		\toprule
		\multicolumn{1}{c}{} & \multicolumn{3}{c}{$10^{-4}$} & \multicolumn{3}{c}{$10^{-6}$} \\
		\cmidrule(lr){2-4} \cmidrule(lr){5-7}
		\LP Method & \# found & gap (\%) & time (s) & \# found & gap (\%) & time (s)\\
		\midrule
		best \PDLP & 801 & 19.39 & 3.70 & 791 & 19.23 & 5.96 \\
		best \IPM  & 794 & 19.25 & 2.19 & 799 & 19.49 & 2.26 \\
		\bottomrule
	\end{tabular}
\end{table}
On \MIPLIB, \PDLP is still slower than \IPM in absolute terms, a result in line with \PDLP literature as most instances in \MIPLIB are small and thus not well suited for \FOMs to outperform \IPMs.
In \cref{tab:pdlp_ipm_wins_losses_size} we count the percentage of times the \IPM or \PDLP variants were faster at finding a solution for all experiments from \cref{tab:appendix:miplib_presets}.
Additionally, we supply the average number of nonzeros of these instances.
\begin{table}[ht!]
    \caption{Percentage of total runs where \LP method found a solution quicker vs average instance nonzeros.}
	\label{tab:pdlp_ipm_wins_losses_size}
	\begin{tabular}{l|cr}
		\toprule
		\LP Method & \% faster & avg. nonzeros\\
		\midrule
		\PDLP & 18.2 & 790\,462.6 \\
		\IPM  & 81.8 & 456\,624.3 \\
		\bottomrule
	\end{tabular}
\end{table}
\PDLP is more likely to outperform on larger or denser models with more non-zeros, a fact that will be reinforced in \cref{sec:numerical_results_esom_pdlp_vs_ipm,sec:numerical_results_ESOM}.

\paragraph{Timing analysis}

Lastly, we present a breakdown of the \FP algorithm by time spent in its major components: presolve, the initial \LP solve (if any), the fix-and-propagate scheme (possibly including repairs), and the final \LP solve.
\Cref{tab:miplib_time_breakdown} shows the average time spent per \FP component for instances where \FP was able to find a solution.
Again we compare all \LP-free presets from \cite{Salvagnin2024_AFixPropagateRepairHeuristicsForMIP} and our additions, split by \LP method.
\begin{table}[ht!]
	\caption{Average percentage breakdown of \FP time by components.}
	\label{tab:miplib_time_breakdown}
	\begin{tabular}{l|ccccc}
		\toprule
		  Methods & reading & presolve & fix-and-propagate & initial \LP & final \LP\\
		\midrule
		  Presets (\LP-free)   & 11.97 & 75.82 & 5.44 & 0.0 & 6.77 \\
		  \LP-based with \PDLP & 1.30 & 9.12 & 1.05 & 87.64 & 0.90 \\
		  \LP-based with \IPM  & 5.86 & 43.63 & 4.25 & 41.71 & 4.55 \\
		\bottomrule
	\end{tabular}
\end{table}
For the \MIPLIB instances, presolve time is relatively significant, especially on the \LP-free methods taking up to 75\% of the total runtime.
For \PDLP-based variants the initial \LP time dominates taking 88\% on average, the \IPM-based methods spend 42\% in the initial \LP.
The fix-and-propagate loop itself is relatively cheap for all variants taking at most 5.5\% of the total runtime.
For \MIPLIB, solving the final \LP is relatively inexpensive taking at most 7\% of total runtime for all aggregations.
We provide a similar breakdown in \cref{sec:numerical_results_ESOM} for the sets miso\_2M and miso\_4M.

\subsection{PDLP vs. IPM performance on large ESOMs}\label{sec:numerical_results_esom_pdlp_vs_ipm}

We compare \PDLP and \IPM on the \LP relaxations of the \ESOMs, crossover disabled.
\Cref{tab:lp_results} shows solution times for \PDLP and \IPM across instance sets and tolerances and we supply the \CPU and \GPU memory used as reported by slurm and ``nvidia-smi''.
\begin{table}[ht!]
\caption{Time to solution (in seconds) and \CPU/\GPU used memory (in \GB) for \COPT 7.1 \PDLP and \IPM with varying tolerances on \ESOMs.}
\label{tab:lp_results}
\begin{tabular}{l|rrrr|rrrr}
\toprule
                 & \multicolumn{2}{c}{{PDLP}} & \multicolumn{2}{c|}{MEM} & \multicolumn{2}{c}{{IPM}} & \multicolumn{2}{c}{MEM}\\
Test Set         & \multicolumn{1}{c}{$10^{-4}$} & \multicolumn{1}{c}{$10^{-6}$} & \CPU & \GPU & \multicolumn{1}{c}{$10^{-4}$} & \multicolumn{1}{c}{$10^{-6}$} & \CPU & \GPU \\
\midrule
miso\_2M    & 13.08  & 22.21   & 0.003 & 1.0 & 29.62    & 31.91 & 1.53 & 0.0\\
miso\_4M    & 9.18   & 26.15   & 0.003 & 1.1 & 72.34    & 86.74 & 2.4 & 0.0\\
miso\_47M   & 104.19 & 166.34  & 9.6 & 3.4 & 1\,002.05  & 1\,216.62 & 18.9 & 0.0 \\
miso\_82M   & 394.54 & 1\,669.78 & 17.2 & 6.4 & 4\,354.58  & 6\,893.95 & 38.3 & 0.0\\
remix\_105M & 147.82 & 552.33  & 21.5 & 6.8 & 10\,067.46 & 12\,893.38 & 113.1 & 0.0 \\
remix\_243M & 436.92 & 2\,131.95 & 49.4 & 13.4 & TIMEOUT  & TIMEOUT & 611.5 & 0.0 \\
\bottomrule
\end{tabular}
\end{table}

All runs used 32 threads, tolerances of $10^{-4}$ and $10^{-6}$, and a time limit of 43{,}200s (12 h).
As instance size increases, \PDLP becomes increasingly competitive.
For the largest models remix\_243M, the \IPM failed to converge within the time limit.
\PDLP benefits more from relaxed accuracy: for large instances, reducing the tolerance to $10^{-4}$ yields speedups of up to a factor of 5, while the \IPM achieves only a factor of 1--2.
\PDLP only stores the system matrix and its required memory scales roughly linearly with matrix size, for our instances \GPU memory stays below 50 \GB.
\IPMs require storing the matrix factor(s) and memory requirements scale closer to quadratically with matrix size reaching more than 600 \GB for our largest models.

The observed advantage of \PDLP stems from both the algorithm and \GPU acceleration.
Both, \PDLP and \IPM are memory bound algorithms.
As \PDLP relies only on vector and sparse matrix-vector operations this effect is even stronger than for the \IPM, where more computations are done relative to memory loaded.
Thus, \PDLP benefits strongly from the high-bandwidth memory of \GPUs.
\GPU implementations of \IPMs quickly exceed the available memory and are more complex to implement efficiently on \GPU.
They still benefit (but less) from the increased memory bandwidth.
Consequently, the relative performance of \PDLP and \IPM depends on available hardware: while \PDLP excels on modern high-bandwidth \GPUs, \IPMs may outperform it on systems with weaker or no \GPU support.

\subsection{Large-scale ESOM experiments}\label{sec:numerical_results_ESOM}

We next evaluated our \FP heuristic on the large-scale \ESOMs using different LP methods and precisions.
We ran all presets and our additions on the miso\_2M and miso\_4M instances (see \cref{tab:appendix:miso2M_presets,tab:appendix:miso2M_lpbased,tab:appendix:miso4M_presets,tab:appendix:miso4M_lpbased}).
In \cref{tab:esom_time_breakdown} we provide a similar time breakdown as in \cref{tab:miplib_time_breakdown} for the sets miso\_2M and miso\_4M.
\begin{table}[ht!]
	\caption{Average percentage breakdown of \FP time by components.}
	\label{tab:esom_time_breakdown}
	\begin{tabular}{ll|ccccc}
		\toprule
		  Set & Methods & reading & presolve & fix-and-propagate & initial \LP & final \LP\\
		\midrule
\multirow{3}{*}{miso\_2M}         & Presets (\LP-free)   & 18.76 & 45.78 & 0.13 & 0.00 & 35.32 \\
                                  & \LP-based with \PDLP & 5.48 & 13.12 & 0.00 & 70.36 & 11.04 \\
                                  & \LP-based with \IPM  & 5.88 & 15.59 & 0.20 & 66.42 & 11.91  \\
		\midrule
\multirow{3}{*}{miso\_4M}         & Presets (\LP-free)   & 9.24 & 22.06 & 0.34 & 0.00 & 68.36 \\
                                  & \LP-based with \PDLP & 3.88 & 9.23 & 0.00 & 54.42 & 32.47 \\
                                  & \LP-based with \IPM  & 1.41 & 3.52 & 0.32 & 83.03 & 11.72 \\
		\bottomrule
	\end{tabular}
\end{table}
The time spend in the initial and final \LP rises for growing instance size taking 67\% and 70\% for \PDLP- and \IPM-based variants respectively on miso\_2M.
For miso\_4M the \IPM based methods even spend 84\% of the time in the initial \LP while for \PDLP this fraction decreases to 54\% as the final \IPM \LP solve becomes more expensive relatively.
Presolve and reading time become less significant for \LP-based variants and the fix-and-propagate logic itself continues to play no role in the overall runtime.

\begin{table}[ht!]
\caption{\ESOM results aggregated by \LP-based and \LP-free \FP strategies.}
\label{tab:miso_lpbased_vs_lpfree}
\begin{tabular}{l|ccc|ccc}
\toprule
                  & \multicolumn{3}{c|}{{miso\_2M}} & \multicolumn{3}{c}{{miso\_4M}} \\
Strategies & \# found & gap (\%) & time (s) & \# found & gap (\%) & time (s) \\
\midrule
Presets            & 100 & 0.05 & 13.6 & 96 & 0.02 & 163.15 \\
Presets (\LP-free) & 100 & 2.91 & 4.30 & 96 & 27.36 & 12.41 \\
\LP based (ours)   & 100 & 0.04 & 11.20 & 96 & 0.02 & 41.58 \\
\midrule
Best                    & 100 & 0.04 & 11.12 & 96 & 0.01 & 77.57 \\
\bottomrule
\end{tabular}
\end{table}

\Cref{tab:miso_lpbased_vs_lpfree} shows \LP-based \FP variants were significantly more successful (while more expensive) on the \ESOMs.
Our additions, as they make use of the fast and low precision \PDLP solutions were also significantly faster than the presets.
We note that the presets performances is mainly driven by the \textit{core} strategy, similar to our \textit{type}, using \IPM with standard accuracy.
We decided for the most successful (compare \cref{tab:appendix:miso2M_lpbased,tab:appendix:miso4M_lpbased}) in terms of gap and instances solved strategy, \textit{type}, combined with the repair mechanism form \cite{Salvagnin2024_AFixPropagateRepairHeuristicsForMIP} and our new branching strategy.
While \textit{type} was not the most successful strategy on \MIPLIB, for our specific instance class it is, displaying the strong variability of primal heuristics performance for diverse sets of \MIPs.
\Cref{tab:uc_results} shows the performance metrics of \textit{type} for all model sets, \LP methods, and precisions with a timelimit of 12 hours.

\begin{table}[ht!]
\caption{Performance of type \FP heuristic on large-scale \ESOMs for \IPM and \PDLP.}
\label{tab:uc_results}
\addtolength{\tabcolsep}{-2.2pt}
\begin{tabular}{ll|ccr|ccr|rr}
    \toprule
    \multicolumn{2}{c}{} & \multicolumn{3}{c}{1e-4} & \multicolumn{3}{c}{1e-6} & \multicolumn{2}{c}{MEM (GB)}\\
    \cmidrule(lr){3-5} \cmidrule(lr){6-8} \cmidrule{9-10}
    Method                 & Testset     & \# found & gap (\%) & time (s) & \# found & gap (\%) & time (s) & CPU & GPU \\
    \midrule
    \multirow{6}{*}{\IPM}  & miso\_2M    & 100 & 0.05 &   12.58 &  99 & 0.05 &   14.10 & 0.003 & 0.0 \\
                           & miso\_4M    & 96 & 0.01 &  73.82 &  96 & 0.01 &  86.54 & 3.0 & 0.0 \\
                           & miso\_47M   & 100 & 0.07 & 1\,391.72 &  100 & 0.11 & 1\,584.05 & 22.1 & 0.0 \\
                           & miso\_82M   &  15 & 0.78 & 7\,707.47 & 16 & 0.65 & 9\,647.04 & 47.3 & 0.0 \\
                           & remix\_105M & 100 & 1.64 & 10\,590.73 & 95 & 1.52 & 13\,252.69 & 58.22 & 0.0 \\
    \midrule
    \multirow{6}{*}{\PDLP} & miso\_2M    & 99 & 0.05 & 6.98     &  100 & 0.05 &   16.17 & 0.003 & 0.6 \\
                           & miso\_4M    & 95 & 0.04 & 18.07&  96 & 0.04 &   34.79 & 0.003 & 0.9 \\
                           & miso\_47M   &  100 & 0.07 & 486.95   &  100 & 0.07 &  558.78 & 20.1 & 3.1 \\
                           & miso\_82M   &  18 & 0.77 & 3\,461.62  &   15 & 0.65 & 5\,281.46 & 42.5 & 6.1 \\
                           & remix\_105M &  96 & 1.74 & 3\,322.69  &   96 & 1.55 & 5\,095.18 & 45.6 & 6.5 \\
                           & remix\_243M &  87 & 1.57 & 12\,587.51 &   82 & 1.48 & 26\,334.98 & 105.8 & 13.1 \\
    \bottomrule
\end{tabular}
\end{table}
%
Most gaps remain below 1\% for the smaller \ESOMs, and below 2\% for two largest sets.
Runtime is dominated by the final barrier solve; the initial LP solve with \PDLP is comparatively cheaper.
As we solve the final \LP using an \IPM, memory requirements for \PDLP based variants also grow, though slightly less than the \IPM counterpart as here the larger initial \LP is also solved by an \IPM.
Memory requirements are smaller than in \cref{tab:lp_results}, as we additionally apply \MIP presolve implemented by Gurobi.
For individual instance sets, either variant can slightly outperform the other, a result which we mostly attribute to performance variability within the method.
E.g., on remix\_105M the \IPM variant produces solutions better by 0.03\% gap and finds four more solutions overall, whereas on miso\_82M, the \PDLP variant found more solutions with 0.01\% smaller gap.
Overall, the heuristic is effective at exploiting low-quality \LP solutions, without sacrificing solution quality but within less time.

\begin{table}[ht!]
\caption{Performance of \GUROBI (1\% target gap) on \ESOMs.}
\label{tab:gurobi_unseen}
\begin{tabular}{lccr}
\toprule
Instance set & \# found & gap (\%) & time (s)  \\
\midrule
miso\_2M    & 100 &  0.40 &     130.81 \\
miso\_4M    &  96 &  0.05 &     881.53 \\
miso\_47M   &  89 &  0.06 & 25\,283.24 \\
miso\_82M   &   0 &     - &          - \\
remix\_105M &   0 &     - &          - \\
remix\_243M &   0 &     - &          - \\
\bottomrule
\end{tabular}
\end{table}

In \cref{tab:gurobi_unseen} we compare our results to \GUROBI (timelimit two days).
Our \PDLP-based \FP heuristic can generate solutions for most instances, often in significantly less time.
\GUROBI, taking a different, global solution approach is slowed down by the initial \LP solve and crossover, failing for the larger sets.
%

\section{Conclusion}

Embedding \FOMs via \PDLP into the \FP heuristic enables systematic exploitation of low-quality \LP solutions without compromising solution quality.
On \MIPLIB 2017 we show that the low quality solutions produced by \FOMs marginally affect solution quality and on the large-scale \ESOMs, this approach achieves gaps below 2\%, often below 1\%, and can solve instances where traditional \IPM-based heuristics or commercial solvers are impractical.
%
%
For future research we believe exploration of more hybrid \LP-based variants to be promising for \FP itself.
The exploitation of \PDLP's cycling convergence behavior to generate multiple \LP reference solutions starting \FP searches as done in \cite{Rothberg_ConcurrentCrossover} for crossover seems promising.
Using \PDLP for the final \LP solve might further speedup the heuristic but solutions would be less accurate.

\paragraph{\footnotesize\normalfont\textbf{Acknowledgements}}{\footnotesize We thank Domenico Salvagnin for making his \FP repair code \cite{Salvagnin2024_AFixPropagateRepairHeuristicsForMIP} available.}
\paragraph{\footnotesize\normalfont\textbf{Funding}}{\footnotesize The described research activities are funded by the Federal Ministry for Economic Affairs and Energy within the project UNSEEN (ID: 03EI1004D). Part of the work for this article has been conducted in the research Campus MODAL funded by the German Federal Ministry of Education and Research (BMBF) (fund numbers 05M14ZAM, 05M20ZBM).}
\paragraph{\footnotesize\normalfont\textbf{Data availability}}{\footnotesize The \ESOMs in this study are available at \cite{ZenodoUNSEEN,ZenodoUNSEENpresolved}, logs and evaluation scripts at \cite{ZenodoFPLogs}. The source code is not publicly available at this time, but compiled binaries can be obtained upon request from kempke@zib.de.}
\paragraph{\footnotesize\normalfont\textbf{Conflict of Interest}}{\footnotesize All authors declare that they have no conflicts of interest.}

\bibliographystyle{mpc_article/spmpsci}      
\bibliography{references_mpc}   

\newpage
\appendix

\section{MIPLIB2017 results}\label{sec:appendix:miplib_res}
\vspace{1.4cm}
\begin{table}[h]
    \rotatebox{90}{%
        \begin{minipage}{0.85\textheight} 
        \caption{Geo. mean of gap (shift 1\%) and time (shift 1s) for \FP variants as reported in \cite{Salvagnin2024_AFixPropagateRepairHeuristicsForMIP} on \MIPLIB.\\Gaps and time are reported for instances where a feasible solution could be found.} \label{tab:appendix:miplib_presets}

\begin{tabular}{l|ccc|ccc|ccc|ccc}
			\toprule\noalign{\smallskip}
			\multicolumn{1}{c}{} & \multicolumn{3}{c}{no repair + old branching} & \multicolumn{3}{c}{repair + old branching} & \multicolumn{3}{c}{no repair + new branching} & \multicolumn{3}{c}{repair + new branching} \\
			\cmidrule(lr){2-4} \cmidrule(lr){5-7} \cmidrule(lr){8-10} \cmidrule(lr){11-13}
			Preset & \# found & gap (\%) & time (s) & \# found & gap (\%) & time (s) & \# found & gap (\%) & time (s) & \# found & gap (\%) & time (s)\\
			\noalign{\smallskip}\midrule\noalign{\smallskip}
						random    & 536 & 63.46 & 1.05 & 654 & 54.75 & 1.09 & 532 & 63.12 & 1.05 & 650 & 54.69 & 1.06 \\
			random2   & 495 & 58.25 & 1.18 & 652 & 55.23 & 1.31 & 486 & 57.22 & 1.24 & 650 & 54.31 & 1.25 \\
			badobj    & 625 & 63.95 & 1.00 & 728 & 57.17 & 1.10 & 625 & 63.95 & 1.02 & 728 & 57.17 & 1.09 \\
\bottomrule
		\end{tabular}
    \end{minipage}%
    }
\end{table}
\newpage
\clearpage
\begin{table}[!ht]
    \setlength{\tabcolsep}{5pt}
    \rotatebox{90}{%
        \begin{minipage}{0.96\textheight} 
            \caption{Geo. mean of gap (shift 1\%) and time (shift 1s) for \LP-based \FP variants on \MIPLIB.\\Gaps and time are reported for instances where a feasible solution could be found.}\label{tab:appendix:miplib_lpbased}
            \begin{tabular}{lll|ccc|ccc|ccc|ccc}
    \toprule\noalign{\smallskip}
    \multicolumn{3}{c}{} & \multicolumn{3}{c}{no repair + old branching} & \multicolumn{3}{c}{repair + old branching} & \multicolumn{3}{c}{no repair + new branching} & \multicolumn{3}{c}{repair + new branching} \\
    \cmidrule(lr){4-6} \cmidrule(lr){7-9} \cmidrule(lr){10-12} \cmidrule(lr){13-15}
    Strategy & method & Tol. & \# found & gap (\%) & time (s) & \# found & gap (\%) & time (s) & \# found & gap (\%) & time (s) & \# found & gap (\%) & time (s)\\
    \noalign{\smallskip}\midrule\noalign{\smallskip}
    \multirow{4}{*}{duals}   & \multirow{2}{*}{\PDLP} & $10^{-4}$ & 516 & 41.33 & 3.34 & 658 & 41.07 & 3.34 & 504 & 35.34 & 3.28 & 659 & 34.74 & 3.27 \\
                             &                        & $10^{-6}$ & 517 & 41.18 & 5.95 & 652 & 39.56 & 5.29 & 505 & 33.95 & 5.86 & 654 & 32.82 & 5.22 \\
\bottomrule
            \end{tabular}
        \end{minipage}%
    }
\end{table}

\newpage

\section{ESOM results}\label{sec:appendix:esoms_res}
\vspace{1.4cm}
\begin{table}[h]
    \rotatebox{90}{%
        \begin{minipage}{0.85\textheight} 
        \caption{Geo. mean of gap (shift 1\%) and time (shift 1s) for \FP variants as reported in \cite{Salvagnin2024_AFixPropagateRepairHeuristicsForMIP} on miso\_2M.\\Gaps and time are reported for instances where a feasible solution could be found.} \label{tab:appendix:miso2M_presets}

        \begin{tabular}{l|ccc|ccc|ccc|ccc}
			\toprule\noalign{\smallskip}
			\multicolumn{1}{c}{} & \multicolumn{3}{c}{no repair + old branching} & \multicolumn{3}{c}{repair + old branching} & \multicolumn{3}{c}{no repair + new branching} & \multicolumn{3}{c}{repair + new branching} \\
			\cmidrule(lr){2-4} \cmidrule(lr){5-7} \cmidrule(lr){8-10} \cmidrule(lr){11-13}
			Preset & \# found & gap (\%) & time (s) & \# found & gap (\%) & time (s) & \# found & gap (\%) & time (s) & \# found & gap (\%) & time (s)\\
			\noalign{\smallskip}\midrule\noalign{\smallskip}
						random    & 100 & 7.34 & 4.34 & 100 & 7.34 & 4.35 & 100 & 5.57 & 4.44 & 100 & 5.65 & 4.48 \\
			random2   & 100 & 6.75 & 4.54 & 100 & 6.75 & 4.52 & 0 & - & - & 0 & - & - \\
			badobj    & 100 & 2.97 & 4.29 & 100 & 2.97 & 4.30 & 100 & 2.97 & 4.31 & 100 & 2.97 & 4.31 \\
\bottomrule
		\end{tabular}
    \end{minipage}%
    }
\end{table}
\newpage
\clearpage
\begin{table}[!ht]
    \setlength{\tabcolsep}{5pt}
    \rotatebox{90}{%
        \begin{minipage}{\textheight} 
            \caption{Geo. mean of gap (shift 1\%) and time (shift 1s) for \LP-based \FP variants on miso\_2M.\\Gaps and time are reported for instances where a feasible solution could be found.}\label{tab:appendix:miso2M_lpbased}
            \begin{tabular}{lll|ccc|ccc|ccc|ccc}
    \toprule\noalign{\smallskip}
    \multicolumn{3}{c}{} & \multicolumn{3}{c}{no repair + old branching} & \multicolumn{3}{c}{repair + old branching} & \multicolumn{3}{c}{no repair + new branching} & \multicolumn{3}{c}{repair + new branching} \\
    \cmidrule(lr){4-6} \cmidrule(lr){7-9} \cmidrule(lr){10-12} \cmidrule(lr){13-15}
    Strategy & method & Tol. & \# found & gap (\%) & time (s) & \# found & gap (\%) & time (s) & \# found & gap (\%) & time (s) & \# found & gap (\%) & time (s)\\
    \noalign{\smallskip}\midrule\noalign{\smallskip}
    \multirow{4}{*}{duals}   & \multirow{2}{*}{\PDLP} & $10^{-4}$ & 100 & 0.31 & 7.23 & 100 & 0.31 & 7.23 & 100 & 0.13 & 7.03 & 100 & 0.13 & 7.05 \\
                             &                        & $10^{-6}$ & 100 & 0.31 & 16.51 & 100 & 0.31 & 16.37 & 100 & 0.13 & 16.06 & 100 & 0.13 & 16.13 \\
\bottomrule
            \end{tabular}
        \end{minipage}%
    }
\end{table}

\newpage
\clearpage

\begin{table}[h]
    \rotatebox{90}{%
        \begin{minipage}{\textheight} 
        \caption{Geo. mean of gap (shift 1\%) and time (shift 1s) for \FP variants as reported in \cite{Salvagnin2024_AFixPropagateRepairHeuristicsForMIP} on miso\_4M.\\Gaps and time are reported for instances where a feasible solution could be found.} \label{tab:appendix:miso4M_presets}
        \begin{tabular}{l|ccc|ccc|ccc|ccc}
			\toprule\noalign{\smallskip}
			\multicolumn{1}{c}{} & \multicolumn{3}{c}{no repair + old branching} & \multicolumn{3}{c}{repair + old branching} & \multicolumn{3}{c}{no repair + new branching} & \multicolumn{3}{c}{repair + new branching} \\
			\cmidrule(lr){2-4} \cmidrule(lr){5-7} \cmidrule(lr){8-10} \cmidrule(lr){11-13}
			Preset & \# found & gap (\%) & time (s) & \# found & gap (\%) & time (s) & \# found & gap (\%) & time (s) & \# found & gap (\%) & time (s)\\
			\noalign{\smallskip}\midrule\noalign{\smallskip}
						random    & 96 & 51.72 & 9.74 & 96 & 51.72 & 9.81 & 96 & 28.34 & 12.09 & 96 & 27.71 & 12.19 \\
			random2   & 96 & 53.03 & 9.30 & 96 & 52.99 & 9.53 & 0 & - & - & 0 & - & - \\
			badobj    & 96 & 53.17 & 11.71 & 96 & 53.17 & 11.75 & 96 & 53.17 & 11.72 & 96 & 53.17 & 11.77 \\
\bottomrule
		\end{tabular}
    \end{minipage}%
    }
\end{table}

\newpage
\clearpage

\begin{table}[!ht]
    \setlength{\tabcolsep}{5pt}
    \rotatebox{90}{%
        \begin{minipage}{0.96\textheight} 
            \caption{Geo. mean of gap (shift 1\%) and time (shift 1s) for \LP-based \FP variants on miso\_4M.\\Gaps and time are reported for instances where a feasible solution could be found.}\label{tab:appendix:miso4M_lpbased}
            \begin{tabular}{lll|ccc|ccc|ccc|ccc}
    \toprule\noalign{\smallskip}
    \multicolumn{3}{c}{} & \multicolumn{3}{c}{no repair + old branching} & \multicolumn{3}{c}{repair + old branching} & \multicolumn{3}{c}{no repair + new branching} & \multicolumn{3}{c}{repair + new branching} \\
    \cmidrule(lr){4-6} \cmidrule(lr){7-9} \cmidrule(lr){10-12} \cmidrule(lr){13-15}
    Strategy & method & Tol. & \# found & gap (\%) & time (s) & \# found & gap (\%) & time (s) & \# found & gap (\%) & time (s) & \# found & gap (\%) & time (s)\\
    \noalign{\smallskip}\midrule\noalign{\smallskip}
    \multirow{4}{*}{duals}   & \multirow{2}{*}{\PDLP} & $10^{-4}$ & 96 & 52.53 & 23.20 & 96 & 52.53 & 23.22 & 96 & 13.45 & 17.85 & 23 & 12.19 & 16.48 \\
                             &                        & $10^{-6}$ & 96 & 52.56 & 40.30 & 96 & 52.56 & 40.36 & 96 & 13.39 & 34.73 & 35 & 12.98 & 30.45 \\
\bottomrule
            \end{tabular}
        \end{minipage}%
    }
\end{table}

\section{cuOpt heuristics-only results}\label{sec:appendix:cuopt}

Results of cuOpt in \textit{heuristics-only}\footnote{\url{https://docs.nvidia.com/cuopt/user-guide/latest/lp-milp-settings.html\#heuristics-only}} mode on \MIPLIB and our \ESOMs.
All instances were presolved using Gurobi 11.0.0 and (presolve times are not included in the cuOpt runtimes) and run on an \ARM Neoverse-V2 \CPU with 72 cores, 480 \GB of \CPU memory, and one \NVIDIA H200 \GPU with 96 \GB of device memory.
Each run executed exclusively using 32 threads; cuOpt's heuristics-only run almost purely on GPU.
Note that cuOpt is not designed for fast-fail solution computation but rather as a standalone heuristic solver.
We enforced tight time-limits on our cuOpt runs to still get some comparability against our primal heuristics.

\begin{table}[h!]
    \caption{Geo. mean of gap (shift 1\%) and time (shift 1s) for cuOpt ``heuristics-only'' on the miso\_2M set with different time limits. Gaps and time are reported for instances where a feasible solution could be found.}\label{tab:appendix:cuopt}
    \begin{tabular}{lc|cc}
        \toprule
        set & timelimit & \# found & gap (\%) \\
        \midrule
        \multirow{2}{*}{\MIPLIB}  & 30s  &  1021 & 11.76 \\
                                  & 300s &  1097 & 6.25  \\
        \midrule
        \multirow{2}{*}{miso\_2M} & 10s &   18 & 2.85 \\
                                  & 30s &  100 & 2.12 \\
        \midrule
        \multirow{2}{*}{miso\_4M} & 30s &   71 & 12.50 \\
                                  & 60s &   91 &  5.29 \\
        \bottomrule
    \end{tabular}
\end{table}

\end{document}